\documentclass[12pt]{amsart} 
\usepackage{amssymb,amsmath,latexsym,enumerate,graphicx,bbm,mathptmx,microtype}
\usepackage{float}
\usepackage{hyperref}
\usepackage{tabularx}
\usepackage{tabularx}
\usepackage{array}

\usepackage[a4paper]{geometry}
\numberwithin{equation}{section}
\newtheorem{theorem}{Theorem}[section]
\newtheorem{corollary}{Corollary}[section]

\newtheorem{exam}{Example}

\newtheorem*{rem}{Remark}

\newcommand\commentout[1]{}

\usepackage{tikz}
\usetikzlibrary{mindmap}

\begin{document}
\date{\today}
\title[\tiny{Common Divisors of Values Polynomials and common factors of indices in a Number Field}]{Common Divisors of Values Polynomials and common factors of indices in a Number Field}
\author{ 
 Mohammed Seddik}
\address{ 
 Mohammed Seddik\\
Universit\'e d'\'Evry Val d'Essonne\\
Laboratoire de MathŽmatiques et Mod\'elisation d'\'Evry (UMR 8071)\\
I.B.G.B.I., 23 Bd. de France, 91037 \'Evry Cedex, France\\}
 
  \email{mohammed.seddik@univ-evry.fr}

\begin{abstract}
Let $\mathbb{K}$ be a number field of degree $n$ over $\mathbb{Q}$. Let $\widehat{\mathbb{A}}$ be the set of integers of $\mathbb{K}$ which are primitive over $\mathbb{Q}$ and $I(\mathbb{K})$ be its index. Gunji and McQuillan defined the following integer  $i(\mathbb{K})=\underset{\theta\in \widehat{\mathbb{A}}}{\text{lcm}}\;i(\theta)$, where $i(\theta)=\underset{x\in \mathbb{Z}}{\text{gcd}}\;F_\theta(x)$ and $F_\theta(x)$ is the characteristic polynomial of $\theta$ over $\mathbb{Q}$. We prove that if $p$ is a prime number less than or equal to $n$ then there exists a number field $\mathbb{K}$ of degree $n$ for which $p$ divides $i(\mathbb{K})$. We compute $i(\mathbb{K})$ for cubic fields. Also we determine $I(\mathbb{K})$ and $i(\mathbb{K})$ for families of simplest number fields of degree less than $7$. We give also answers to questions one and two in  \cite{Kihel}. 
\end{abstract}




\maketitle

\section{\textbf {Introduction and statement of main results}}

Let $\mathbb{K}$ be a number field of degree $n$ over $\mathbb{Q}$ and let $\mathbb{A}$ be its ring of integers. Denote by $\widehat{\mathbb{A}}$ the set of primitive elements of $\mathbb{A}$. For any $\theta\in \mathbb{A}$ we denote $F_\theta(x)$ the characteristic polynomial of $\theta$ over $\mathbb{Q}$. Let $D_\mathbb{K}$ be the discriminant of $\mathbb{K}$. It is well known that if $\theta\in \widehat{\mathbb{A}}$, the discriminant of $F_\theta(x)$ has the form 
\begin{equation}\label{equ1}
D(\theta)=I(\theta)^2D_\mathbb{K}
\end{equation}

 where $I(\theta)=[\mathbb{A}:\mathbb{Z}[\theta]]$ is called the index of $\theta$. Let
 
 \begin{equation}\label{equ2}
 I(\mathbb{K})=\underset{\theta\in \widehat{\mathbb{A}}}{\text{gcd}}\;I(\theta).
\end{equation} 
 A prime number $p$ is called a common divisor of indices in $\mathbb{A}$ or sometimes a common index divisor, if $p\mid I(\mathbb{K})$. Dedekind was the first one to show the existence of common divisor of indices. He exhibited an example of a number field of degree $3$ in which $2$ is a common divisor of indices \cite[pp.183-184]{Alaca}. Bauer \cite{Bauer} showed that if $p<n$ then there exists a number field of degree $n$ in which $p$ is a common index divisor. Zylinski \cite{Zylinski} showed the necessity of this condition, if $p$ is  common index divisor then $p<n$. Hensel \cite{Hensel} has given a necessary and sufficient condition on a prime $p$ to be a common divisor of indices in a number field $\mathbb{K}$. This condition depends upon the splitting of the prime $p$ in $\mathbb{K}$, which make Hensel's Theorem not easy to apply in general.\\
 \vspace{0.2cm}
 
Let $\theta\in\mathbb{A}$ and $i(\theta)=\underset{x\in \mathbb{Z}}{\text{gcd}}\;F_\theta(x)$. Gunji and McQuillan \cite{Gunji} defined the following integer 
\begin{equation}\label{equ3}
i(\mathbb{K})=\underset{\theta\in \widehat{\mathbb{A}}}{\text{lcm}}\;i(\theta),
\end{equation}
and they showed that for $m$ square free rational integer,
\begin{center}
$i(\mathbb{Q}(\sqrt{m}))=
\begin{cases}
   2 & \text{if $m \equiv 1 \;mod\; 8$}\\
   1 & \text{otherwise.}   
\end{cases}$
 \end{center}

Mac Cluer \cite{Mac Cluer} showed that $i(\mathbb{K})>1$ if and only if there exists a prime number $p\leq n$ having at least $p$ distinct prime ideal factors in $\mathbb{A}$, each of these primes and only these primes are divisors of $i(\mathbb{K})$. Ayad and Kihel \cite {Kihel} gives two important theorems. It is shown that there exist a primitive integer $\theta$ called a good element such that $i(\mathbb{K})=i(\theta)$ and an algorithm is given for the computation of such an integer and showed that if $p$ is a common factor divisor then $p\mid i(\mathbb{K})$. The converse is shown to be false in general. However, the following result is proved. Suppose that $\mathbb{K}$ is a Galois extension over $\mathbb{Q}$. Let $1\leq d<n$ be the greatest divisor of $n$ and let $p>d$ be a prime number, $p\neq n$. Then $p$ is a common index divisor if and only if $p\mid i(\mathbb{K})$. As a consequence, we obtain that if $\mathbb{K}/\mathbb{Q}$ is cyclic of prime degree $l$, then $p\neq l$ is a common index divisor if and only if $p\mid i(\mathbb{K})$. Let $p$ be a prime number. Let $v_p(i(\mathbb{K}))$ be the valuation of $p$ in $i(\mathbb{K})$.
\vspace{0.2cm}

In 1926, Ore \cite{Ore} conjectured that p-adic valuation $v_p(I(\mathbb{K}))$ is not determined only by the splitting type of $p$ in $\mathbb{A}$. Engstrom \cite{Engstrom} proved that if $n\leq 7$, then the splitting type determines the p-adic valuation $v_p(I(\mathbb{K}))$. He gave examples of number fields $\mathbb{K}_1$ and $\mathbb{K}_2$ of degree 8 in which the prime $3$ has the same splitting type, but $v_3(I(\mathbb{K}_1))\neq v_3(I(\mathbb{K}_2))$. \'Sliwa \cite{Sliwa} proved that if $p$ is not ramified, then $v_p(I(\mathbb{K}))$ is determined by the splitting type of $p$ in $\mathbb{K}$. Similar to this conjecture of Ore, Ayad and Kihel \cite{Kihel} ask the following question,
\vspace{0.2cm}

Suppose that the splitting of $p$ in $\mathbb{K}$ as a product of prime ideals is given by $p\mathbb{A}=P_1^{e_1}\cdots P_r^{e_r}$, where $r\geq p$. Let $f_i$ be the inertial degree of the ideal $P_i$, for $i=1,\ldots,r$. Can one compute $v_p(I(\mathbb{K}))$ in terms of integers $r,\,e_i$ and $f_i$. In other words, is $v_p(i(\mathbb{K}))$ completely determined by splitting type of $p$ ?. We give answer of this question, we gives examples of number fields $\mathbb{K}_1$ and $\mathbb{K}_2$ of degree $6$ in which the prime $3$ has the same splitting type $P_1P_2$ but $v_3(I(\mathbb{K}_1))\neq v_3(I(\mathbb{K}_2))$.
\vspace{0.2cm}

 They also showed that, if $\mathbb{K}_1$ and $\mathbb{K}$ be number fields such that $\mathbb{K}_1\subseteq \mathbb{K}$ and let $m=[\mathbb{K}:\mathbb{K}_1]$, then $mv_p(i(\mathbb{K}_1))\leq v_p(i(\mathbb{K}))$. Answer the question 2 in \cite{Kihel}, we show that the following statements are not equivalent,
\begin{enumerate}
\item $mv_p(i(\mathbb{K}_1))= v_p(i(\mathbb{K}))$.
\item For any integer $\beta$ of $\mathbb{K}$, if $v_p(i(\beta))=v_p(i(\mathbb{K}))$, then there exists an integer $\alpha$ of $\mathbb{K}_1$ such that $\beta\equiv \alpha \;mod\;p$.
\end{enumerate}
\vspace{0.2cm}

We now state the main result of this paper.

\begin{theorem}\label{thm1}
Let $\mathbb{K}$ be a number field of degree $n$ over $\mathbb{Q}$ and $p$ a prime number.

 If $p\leq n$ then, there exists a number field $\mathbb{K}$ of degree $n$ in which $p\mid i(\mathbb{K})$.

\end{theorem}
Now we will calculate $i(\mathbb{K})$ for cubic number fields. Let $\mathbb{K}$ be a cubic field. We can suppose that $\mathbb{K}=\mathbb{Q}(\theta)$, where $\theta$ is a root of an irreducible polynomial of the type 
$$f(x)=x^3-ax+b,\;a,b\in  \mathbb{Z}.$$
The discriminant of $f(x)$ is $\Delta=4a^3-27b^2$. If for any prime number $p$ we have 
$$v_p(a)\geq 2 \quad\text{and}\quad v_p(b)\geq 3 $$
then $\theta/p$ is an algebraic integer whose equation is $x^3-(a/p^2)x+b/p^3=0$. Therefore, we can assume that for any prime number $p$,
$$v_p(a)< 2 \quad\text{or}\quad v_p(b)< 3 .$$
Let, $$ s_3=v_3(\Delta),\quad \Delta_3=\Delta/3^{s_3}.$$
It is well known that if $\mathbb{K}$ is a cubic field , then $I(\mathbb{K})=1$ or $2$ and 
$$I(\mathbb{K})=2\Longleftrightarrow a\;\text{odd}\;,b\;\text{even},\;s_2\;\text{even}\;\text{and}\;\Delta_2\equiv 1\;mod\;8.$$
see \cite[Theorem 4]{Llorente}.
\begin{theorem}\label{thm5}
Let $\mathbb{K}=\mathbb{Q}(\theta),\;\theta^3-a\theta+b=0,$ be a cubic field. Then the common divisors of values polynomials is given by $$i(\mathbb{K})=2^\alpha3^\beta,$$
 where
 \[
\alpha=
         \begin{cases}
           1 & \text{if}\; 1=v_2(a)<v_2(b) \;\text{or}\; a\not\equiv b\;mod\;2,\\
           0 & \text{else}.   
\end{cases}
\]
 \[
\beta=
\begin{cases}
        1 & \text{ if} \;(a\equiv 3\;mod\;9,\; b^2\equiv a+1\;mod\;27,\; s_3 > 6 \; \text{even},\; \Delta_3\equiv 1\;mod\;3)\\
          & \text{or}\; (a\equiv 1\;mod\;3,\; 3\mid b),\\
        0 & \text{else}.   
\end{cases}
\]

\end{theorem}

As a particular case of the above, we get $i(\mathbb{K})$ for $\mathbb{K}$ is a pure cubic field.

\begin{corollary}\label{coro6}
Let $\mathbb{Q}(\sqrt[3]{d})$ be an pure cubic field, then
\begin{center}
$i(\mathbb{Q}(\sqrt[3]{d}))=
\begin{cases}
   2 & \text{if \; d odd},\\
   1 & \text{if \; d even}.   
\end{cases}$
 \end{center}
\end{corollary}
Let us consider the family of cyclic cubic fields $\mathbb{L}_m$ generated by a root of the polynomial 
$$l_m(x)=x^3-mx^2-(m+3)x-1.$$
This family were discussed by Shanks \cite{shanks} and Jager \cite[p 63-73]{Jager}.

\begin{theorem}\label{thm6}
Let a simplest cubic fields given by $\mathbb{L}_m$ generated by a root of the polynomial $l_m$. Then $$I(\mathbb{L}_m)=1.$$
\[
i(\mathbb{L}_m)=
\begin{cases}
   3 & \text{if}\quad m\equiv39,120,201\;mod\;243,\\
   1 & \text{otherwise}.   
\end{cases}
 \]

\end{theorem}
For $m\in\mathbb{Z}$ with $m\not\in\{0,\pm 3\}$ let $P_m(x)=x^4-mx^3-6x^2+mx+1$. Let $\theta$ be a root of $P_m(x)$, then each field in infinite parametric family of number fields $\mathbb{K}_m=\mathbb{Q}(\theta)$ is called a simplest quartic field. If $m\not\in\{0,\pm 3\}$ then $P_m(x)$ is irreducible over $\mathbb{Q}$ and defined a totally real cyclic number field of degree $4$, see \cite[Proposition 6]{Gras}. Note that $\mathbb{K}_m=\mathbb{Q}(\theta)=\mathbb{K}_m=\mathbb{Q}(-\theta)$, that is we can assume that $m>0$ and $m\neq 3$.
\begin{theorem}\label{thm7}
Let $m>0$ and $m\neq 3$. Suppose the $m^2+16$ is not divisible by an odd square. Consider a simplest quartic fields  generated by a root of the polynomial $P_m$. Then 

\[
I(\mathbb{K}_m)=
\begin{cases}
   2 & \text{if \; m odd},\\
   1 & \text{if \; m even}.   
\end{cases}
 \]
 
  \[
   i(\mathbb{K}_m)=
   \begin{cases}
      1 & \text{if}\; 1\leq v_2(m)\leq 3,\\
      4 &  \text{otherwise}. 
        
   \end{cases}
    \]
\end{theorem}

Let us consider the family of sextic field $\mathbb{H}_m$ generated by a root of the polynomial 
$h_m(x)=x^5+m^2x^4-(2m^3+6m^2+10m+10)x^3+(m^4+5m^3+11m^2+15m+5)x^2+(m^3+4m^2+10m+10)x+1.$
 This polynomial were discussed by Emma Lehmer \cite{Lehmer}. The fields $\mathbb{H}_m$ were also investigated by R. Schoof and L. Washington \cite{Washington} and H. Darmon \cite{ Darmon} for prime conductors $m^4+5m^3+15m^2+25m+25$.

\begin{theorem}\label{thm8}(Lehmer's quintics)\\
Let $m\in\mathbb{Z}$ and suppose that $p^2\nmid m^4+5m^3+15m^2+25m+25$ for any prime $p\neq 5$. Consider a simplest quintic fields given by $\mathbb{H}_m$ generated by a root of the polynomial 
$h_m(x)$. Then $$I(\mathbb{H}_m)=1,$$

\begin{center}
   $i(\mathbb{H}_m)=
   \begin{cases}
      5 & \text{if \; $m\equiv 2\;mod\;5$},\\
      1 & \text{otherwise}. 
        
   \end{cases}$
    \end{center} 
\end{theorem}

Assume $m\not\in \{-8,-5,-3,0\}$. Let us consider the family of sextic cyclic fields $\mathbb{S}_m$ generated by a root of the polynomial 
$$s_m(x)=x^6-2mx^5-(5m+15)x^4-20x^4-20x^3+5mx^2+(2m+6)x+1.$$
This family of fields is called the simplest sextic fields, having a couple of nice
properties, detailed in G. Lettl, A. Petho and P. Voutier \cite{Lettl}.

\begin{theorem}\label{thm10}
Let $\mathbb{S}_m$ be a simplest sextic field. Then
\[I(\mathbb{S}_m)=1,\]
 \[i(\mathbb{S}_m)=2^\alpha3^\beta,\]
  where
 \[
 \alpha=
          \begin{cases}
            3 \;\text{or}\;4& \text{if}\;m\equiv0,5\;mod\;8,\;3\nmid m  \\
            &\text{or} \; m\equiv0,21\;mod\;24,\\
            0 & \text{else}.   
 \end{cases}
 \]
 
 \[ 
 \beta=
 \begin{cases}
         2 & \text{ if} \;m\equiv39,120,201\;mod\;243,\\
         
         0 & \text{else}.   
 \end{cases}
 \]
    
\end{theorem}


\begin{thebibliography}{2}
\bibitem{Kihel}
M. \textsc{Ayad} and O. \textsc{Kihel} \emph{Common divisors of values of polynomials and common factors of indices in a number field}, Int. J. Number Theory {\bf 7} (2011), no. 5, 11731194.

\bibitem{Bauer}
N. \textsc{Bauer}. \emph{Uber den ausserwesentlicher discriminamtenteiler algebraischer korper}, Math. Ann. {\bf 64} (1907), 573.

\bibitem{Zylinski}
E. \textsc{Zylinski}. \emph{Zur Theorie der ausserwesentlicher discriminamtenteiler algebraischer korper }, Math. Ann. {\bf 73} (1913) 273-274.

\bibitem{Gunji}
H. \textsc{Gunji}.D. L. \textsc{McQuillan} \emph{On a class of ideals in an algebraic number field}, J. Number Theory {\bf 2} (1970), 207-222.

\bibitem{Mac Cluer}
C. R. \textsc{MacCluer}. \emph{Common divisors of values of polynomials}, J. Number Theory {\bf 3} (1971), 33-34.

\bibitem{Engstrom}
H. T \textsc{Engstrom}. \emph{On the common index divisors of an algebraic field}. Trans. Amer. Math. Soc. {\bf 32} (1930). 223-237.

\bibitem{Sliwa}
J. \textsc{\'Sliwa}. \emph{On the essential discriminant divisor of an algebraic number field}. Acta Arith. {\bf 42} (1983). 57-72.

\bibitem{Llorente}
P. \textsc{Llorente}, \& E. \textsc{Nart}. \emph{Effective determination of the rational primes in a cubic field}, Proc. Amer. Math. Soc. {\bf 87} (1983), 579-585.

\bibitem{Cassels}
J. W. S \textsc{Cassels}. \emph{Local fields}. Cambridge University Press, 1986.

\bibitem{Hass}
H. \textsc{Hass}. \emph{Arithmetische Theorie der kubischen Zahlkörper auf klassenkörpertheoretischer Grundlage}. Math. Zeit. {\bf 31} (1930), 565-582.


\bibitem{Nar}
W. \textsc{Narkeiwicz}. \emph{Elementary and Analytic Theory of Algebraic Numbers }. Second edition. Springer-Verlag, Berlin, 1990.

\bibitem{Delone}
B. N \textsc{Delone} and D. K \textsc{Faddeev} \emph{The theory of irrationalities of the third degree }. Trans. Math. Monographs. vol 10. Amer. Math. Soc. Providence. R. I. 1964.

\bibitem{Alaca}
S. \textsc{Alaca}, K.S. \textsc{Williams}, \emph{Introductory Algebraic Number Theory}, Cambridge University Press, 2004.

\bibitem{shanks}
D. \textsc{Shanks}. \emph{The simplest cubic fields}. Mathematics of computation, {\bf 28} (1974), 1137-1152.

\bibitem{Nagell}
T. \textsc{Nagell}. \emph{Quelques résultat sur les diviseurs fixes de l'index des nombres entiers d'un corps algébrique}. Ark. Mat. {\bf 6} (1965) 269-289.

\bibitem{Hensel}
K. \textsc{Hensel}. \emph{ Theorie der algebraischen zahlen}. Leipzig 1908.

\bibitem{Williams}
K. S \textsc{Williams}. \emph{ Integers of biquadratic fields }. Canad. Math. Bull. Vol. {\bf 13} (4), 1970.

\bibitem{Mayer}
D. C. \textsc{Mayer}, \emph{Multiplicities of dihedral discriminant}, Math. Comp. 58 (1992) 831-847.

\bibitem{Ore}
O. \textsc{Ore}, \emph{Über den Zusammenhang zwischen den definierenden Gleichungen und der Idealtheorie in algebraischen Körpern}, Math. Ann. {\bf 69} (1926) 313-352.

\bibitem{Kim}
H. K. \textsc{Kim} and J. S. \textsc{Kim}, \emph{Computation of the different of the simplest quartic fields}, Manuscript, 2003.

\bibitem{Gras}
M. N. \textsc{Gras}, \emph{Table numérique du nombre de classe et des unités des extensions cycliques réelles de degré 4 de $\mathbb{Q}$}, Publ. Math. Fac. Besançon, fasc 2, 1977-1978.

\bibitem{Gaal}
I. \textsc{Gaal} and M. \textsc{Pohst}, \emph{Power integral bases in a parametric family of
totalty real quintics}, Math. Comput, {\bf 66} (1997), 1689-1696.

\bibitem{Funakura}
T. \textsc{Funakura}, \emph{On integral bases of pure quartic fields}, Math. J. Okayama Univ. {\bf 26} (1984), 27-41.

\bibitem{Lehmer}
E. \textsc{Lehmer}, \emph{Connection between Gaussian periods and cyclic units,} Math. Comput, {\bf 50} (1988), 535-541.


\bibitem{Washington}
R. \textsc{Schoof} and L. \textsc{Washington}, \emph{Quintic polynomials and real cyclotomic
fields with large class numbers,}Math. Comput, {\bf 50} (1988), 543-556.

\bibitem{Darmon}
H. \textsc{Darmon}, \emph{Note on a polynomial of Emma Lehmer,}Math. Comput,
{\bf 56} (1991), 795-800.

\bibitem{Lettl}
G. \textsc{Lettl}, A. \textsc{Petho}, P. \textsc{Voutier} \emph{On the arithmetic of simplest sextic fields and related Thue equations,}in "Number Theory", ed. by K. Gyory, A Petho, y.T. S6s, Walter de Gruyter, Berlin-New York, 1998, pp. 331-348.

\bibitem{Jager}
H. \textsc{Jager}, \emph{Number Theory Noordwijkerhout}, Springer Verlag, 1984.



\end{thebibliography}
\end{document}